\documentclass[12pt,leqno]{article}
\usepackage{amsmath,amsthm,graphicx,color}

\def\init{\setcounter{equation}{0}}
\newtheorem{theorem}{Theorem}[section]

\newcommand{\R}{{\bf R}}

\newtheorem{lemma}{Lemma}[section]

\newcommand{\e}{{\varepsilon}}

\newlength{\figboxwidth}             
\newcommand{\makefig}[3]{
        \begin{figure}[htb]
        \refstepcounter{figure}
        \label{#2}
        \begin{center}
                #3~\\
                \smallskip
                Figure \thefigure.  #1
        \end{center}
        \medskip
        \end{figure}
}

\title{A new approach to  hyperbolic inverse problems II (Global step)				
\author{G.Eskin, \ \ \  Department of Mathematics, UCLA,\\ Los Angeles,
CA 90095-1555, USA. \ E-mail: eskin@math.ucla.edu}
}

\begin{document}

\maketitle
\begin{abstract}
We study the inverse problem for the second order self-adjoint hyperbolic 
equation with the boundary data given on a part of the boundary.  This paper is
the continuation  of the author's paper [E].

In [E]  we presented the crucial local step of the proof.  In this paper we
prove the global step.   Our method is  a modification of the BC-method  with some 
new  ideas.   In particular,  the way  of the determination of the metric is new.  
\end{abstract}

\section{Introduction.}
\label{section 1}
\init

Let $\Omega$  be a bounded domain in $\R^n,\ n\geq 2,$   with
smooth boundary $\partial\Omega$.  Consider the hyperbolic equation of the form:
\begin{eqnarray}                               \label{eq:1.1}
Lu\stackrel{def}{=}
\frac{\partial^2 u}{\partial t^2}
+\sum_{j,k=1}^n\frac{1}{\sqrt{g(x)}}
\left(-i\frac{\partial}{\partial x_j}+A_j(x)
\right)
 \sqrt{g(x)}g^{jk}(x)
\left(-i\frac{\partial}{\partial x_k}+A_k(x)\right)u
\nonumber
\\
+V(x)u=0
\end{eqnarray}
in $\Omega\times(0,T_0)$  with 
$C^\infty(\overline{\Omega})$ coefficients.
Here  
$ \|g^{jk}(x)\|^{-1}$ is the 
metric tensor,
 $g(x)=\det\|g^{jk}\|^{-1}$.  We assume that
\begin{equation}                                 \label{eq:1.2}
u(x,0)=u_t(x,0)=0 \ \ \mbox{in}\ \ \ \Omega,\ \ \
\end{equation}
\begin{equation}
u\left|_{\partial\Omega\times(0,T_0)}\right. = f(x,t).
\end{equation}
Denote by $\Lambda$ the Dirichlet-to-Neumann (D-to-N)   operator,
i.e.
\begin{equation}                                 \label{eq:1.4}
\Lambda f=\sum_{j,k=1}^n g^{jk}(x)\left(\frac{\partial u}{\partial x_j}+
iA_j(x)u\right)\nu_k
\left(\sum_{p,r=1}^n g^{pr}(x)\nu_p\nu_r\right)^{-\frac{1}{2}}
{\Huge |}_{\partial\Omega\times(0,T_0)},
\end{equation}
where
$\nu=(\nu_1,...,\nu_n)$  is the unit exterior normal to $\partial\Omega$
with respect to the Euclidean metric.  Let  $\Gamma_0$  be an open subset of
$\partial\Omega$.  We say that the D-to-N operator is given on
$\Gamma_0\times(0,T_0)$  if $\Lambda f{\Huge |}_{\Gamma_0\times(0,T_0)}$  is
known for all smooth $f(x,t)$ with supports in $\Gamma_0\times (0,T_0]$.

Denote by $G_0(\overline{\Omega})$  the group of all complex-valued 
functions $c(x)$ such that $c(x)\neq 0$ in $\overline{\Omega}$
and $c(x)=1$ on $\overline{\Gamma_0}$.  We say that potentials
$A(x)=(A_1(x),...,A_n(x))$  and $A'(x)=(A_1'(x),...,A_n'(x))$ are
gauge equivalent if there exists  $c(x)\in G_0(\overline{\Omega})$
such that
\[
A_j'(x)=A_j(x)-ic^{-1}(x)\frac{\partial c}{\partial x_j},\ \ \ \
1\leq j\leq n.
\]
Note that if $Lu=0$  then
\begin{equation}                             \label{eq:1.5}
u'=c^{-1}(x)u
\end{equation}
satisfies  the equation $L'u'=0$  where $L'$ has the form (\ref{eq:1.1})
with $A_j(x)$ replaced by $A_j'(x), \  1\leq j\leq n$.  We shall call 
(\ref{eq:1.5})  the gauge transformation.

We shall prove the following theorem:
\begin{theorem}                    \label{theo:1.1} 
Let $L^{(p)},p=1,2$,  be two operators of the form  (\ref{eq:1.1}) in
domains $\Omega^{(p)},p=1,2,$  respectively.  
Let $\Gamma_0\subset\partial\Omega^{(1)}\cap\partial\Omega^{(2)}$ and
let $\Lambda^{(p)},p=1,2,$
be the D-to-N operators corresponding to $L^{(p)}, p=1,2.$  Assume that 
$L^{(p)}$
are self-adjoint,   i.e. coefficients $A_1^{(1)}(x),...,A_n^{(1)}(x),V^{(1)}(x)$
and 
$A_1^{(2)}(x),...,A_n^{(2)}(x),$
\\
$V^{(2)}(x)$
are real-valued.  Suppose $T_0>2\max_{x\in \bar{\Omega}^{(1)}}d_1(x,\Gamma_0)$,
where $d_1(x,\Gamma_0)$  is the distance in $\overline{\Omega_1}$  with
respect to the metric $\|g_1^{jk}(x)\|^{-1}$ from $x\in \overline{\Omega^{(1)}}$
to $\Gamma_0$.
Suppose that the D-to-N operators $\Lambda^{(1)}$ and $\Lambda^{(2)}$  are equal
on $\Gamma_0\times(0,T_0)$ for all $f$ with $\mbox{supp\ } f\subset
\Gamma_0\times (0,T_0]$.  Then there exists a diffeomorphism $\varphi$
of $\overline{\Omega_2}$ onto $\overline{\Omega_1}$, 
$\varphi=I$ on $\Gamma_0$,  and there exists a gauge transformation 
$c(x)\in G_0(\overline{\Omega^{(1)}}$  such that 
$c\circ\varphi\circ L^{(2)}=L^{(1)}$ in $\Omega^{(1)}$.
\end{theorem}  

An important case of the inverse problems with boundary data on a part of
the boundary are the inverse problems in domains  with obstacles.
In this case $\Omega=\Omega_0\setminus\cup_{r=1}^m \overline{\Omega_r}$,  where
$\Omega_0$ is diffeomorphic to a ball,  
$\overline{\Omega_1},...,\overline{\Omega_m}$
are smooth nonintersecting domains in $\Omega_0$  called obstacles,
$\Gamma_0=\partial\Omega_0$ and  zero Dirichlet boundary conditions 
hold on $\partial\Omega_r,\ 1\leq r\leq m$ (c.f. [E1]).

The first result on the inverse problems with the data on a part
 of the boundary was obtained in [I].
The general self-adjoint case was studied by the BC-method (see [B1], [B2],
[K], [KK], [KKL], [KL1]).  The present paper is a continuation of the paper [E]
(see also [E2]).
In [E] the crucial local step was considered,  i.e. the unique determination
of the coefficients of (\ref{eq:1.1}) modulo  a diffeomorphism and a gauge
transformation near $\Gamma_0$.

In this paper we complete the proof of Theorem \ref{theo:1.1}.  In \S 2 we state 
the main results proven in [E] and prove the extension lemma.
In \S 3 we refine the results of \S 2,   
and in \S 4 we complete the proof of Theorem \ref{theo:1.1}.

\section{The summary of the local step and the extension lemma.}
\label{section 2}
\init

Let $L^{(p)}, p=1,2,$  be two operators of 
the form (\ref{eq:1.1}) in $\Omega^{(p)}\times(0,T_0),\ L^{(p)}u_p=0$ in
$\Omega^{(p)}\times (0,T_0),\ u_p(x,0)=u_{pt}(x,0)=0,\ 
x\in \Omega^{(p)},\ u_p|_{\Gamma_0\times(0,T_0)}=f,\ ,p=1,2.$

Let $\Gamma$ be an open connected subset of $\Gamma_0$ and let
$x=(x',x_n)$  be a system of coordinates in a neighborhood $V\subset \R^n$
of $\Gamma$ such that $x_n=0$ is the equation 
of $\Gamma$ and $ x'=(x_1,...,x_{n-1})$ are local coordinates on $\Gamma$.  
Introduce semigeodesic coordinates in $V$ corresponding to $\Gamma$ and to
 the metric
$\|g_p^{jk}\|^{-1},p=1,2$:
\begin{equation}                                \label{eq:2.1}
y=\varphi_p(x).
\end{equation}
Note that $\varphi_p(x)=(\varphi_{p1}(x),...,\varphi_{pn}(x)),p=1,2,$
satisfy  the following differential equations (see [E],  page 817)
\begin{eqnarray}                                    \label{eq:2.2}
\sum_{j,k=1}^n g_p^{jk}(x)\frac{\partial\varphi_{pn}}{\partial x_j}
\frac{\partial\varphi_{pn}}{\partial x_k}=1,\ \ 0\leq x_n<\delta,
\\
\varphi_{pn}(x',0)=0,
\nonumber
\end{eqnarray}

\begin{eqnarray}                                    \label{eq:2.3}
\sum_{j,k=1}^n g_p^{jk}(x)\frac{\partial\varphi_{pn}}{\partial x_j}
\frac{\partial\varphi_{pr}}{\partial x_k}=0,
\\
\varphi_{pr}(x',0)=x_r, \ \ 1\leq r\leq n-1.
\nonumber
\end{eqnarray}
Denote by $\hat{L}^{(p})$ the operator $L^{(p)}$ in semigeodesic
coordinates:
\begin{eqnarray}                    
\hat{L}^{(p)}v_p\stackrel{def}{=}
\frac{\partial^2 v_p(y,t)}{\partial t^2}
+\frac{1}{\sqrt{\hat{g}_p(y)}}\left(-i\frac{\partial}{\partial y_n}+\hat{A}_{pn}(y)\right)
\sqrt{\hat{g}_p(y)}
\left(-i\frac{\partial }{\partial y_n}+\hat{A}_{pn}\right)v_p
\nonumber
\\
+\sum_{j,k=1}^{n-1}\frac{1}{\sqrt{\hat{g}_p(y)}}
\left(-i\frac{\partial}{\partial y_j}+\hat{A}_{pj}(y)\right)
\sqrt{\hat{g}_p}\hat{g}_p^{jk}(y)
\left(-i\frac{\partial}{\partial y_k}+\hat{A}_{pk}(y)\right)v_p(y,t)
\nonumber
\\ 
+\hat{V}_p(y)v_p(y,t)=0,
\nonumber
\end{eqnarray}
where $v_p(\varphi_p(x),t)=u_p(x,t),\ L^{(p)}u_p=0,\ p=1,2$.
Denote
\begin{equation}                         \label{eq:2.4}
u_p^{(1)}(y,t)=(\hat{g}_p(y',y_n))^{\frac{1}{4}}e^{-i\psi_p(y',y_n)}
u_p(\varphi_p^{-1}(y),t),\ \ \ p=1,2,
\end{equation}
where the gauge transformation $c_p(y)= e^{-i\psi_p(y)}\in G_0(\overline{U})$
is such that 
\begin{equation}                          \label{eq:2.5}
A_{pn}^{(1)}(y)=\hat{A}_{pn}(y)+\frac{\partial \psi_p}{\partial y_n}=0,\ \ \ \
\psi_p(y',0)=0.
\end{equation}
Then $u_p^{(1)}$ satisfies the equation (c.f. [E],  page  819):
\begin{eqnarray}                               \label{eq:2.6}
L_1^{(p)}u_p^{(1)}\stackrel{def}{=}\frac{\partial^2 u_p^{(1)}}{\partial t^2}-
\frac{\partial^2 u_p^{(1)}}{\partial y_n^2}
\ \ \ \ \ \ \ \ \ \ \ \ \ \ \ \ \ \ \  
\\
 +
\sum_{j,k=1}^{n}\left(-i\frac{\partial}{\partial y_j}+A_{pj}^{(1)}(y)\right)
\hat{g}_p^{jk}(y) \left(-i\frac{\partial}{\partial y_k}+A_{pk}^{(1)}(y)\right)u_p^{(1)}
\nonumber
\\
+V_{p1}u_p^{(1)}=0,\ \ \ \ p=1,2.
\nonumber
\end{eqnarray}
Let $T\in (0,T_0)$  be  small.  As in [E] denote by $\mathcal{D}_p(\Gamma\times(0,T))$
the forward domain of influence of $\overline{\Gamma}\times[0,T]$  for 
$y_n\geq 0,p=1,2,$  and let 
$\overline{\mathcal{G}_p}
=\{y':(y',y_n,t)\in \mathcal{D}_p(\Gamma\times(0,T)),y_n=0,t=T\}$.
It follows from the proof of Lemma 2.4 in [E]  that $\mathcal{G}_1=
\mathcal{G}_2\stackrel{def}{=}\mathcal{G}$ since $\Lambda^{(1)}=\Lambda^{(2)}$
on $\Gamma_0\times(0,T)$.

Analogously,  let $\mathcal{D}_p(\mathcal{G}\times(0,T))$ be the forward domain of
influence of $\overline{\mathcal{G}}\times[0,T]$ for $y_n\geq 0,p=1,2$.
Let $Y_{20}^{(p)}$ be the intersection of $\mathcal{D}_p(\mathcal{G}\times(0,T))$
with the plane $T-y_n-t=0$  and let $X_{20}^{(p)}$ be the part
of $\mathcal{D}_p(\mathcal{G}\times(0,T))$ below $Y_{20}^{(p)},p=1,2$ 
(c.f. [E],  page 819).
Note that $\Lambda^{(1)}=\Lambda^{(2)}$ on $\Gamma_0\times(0,T_0)$  implies that 
$\hat{\Lambda}^{(1)}=\hat{\Lambda}^{(2)}$ on $\mathcal{G}\times(0,T)$,  where 
$\hat{\Lambda}^{(p)}f'=\frac{\partial u_{p}^{(1)}}{\partial y_n}{\Huge |}_{y_n=0}$
is the D-to-N operator corresponding to $L_1^{(p)}$ (see [E], page 819).
Here $f'=u_p^{(1)}{\Huge |}_{y_n=0},u_p^{(1)}=u_{pt}^{(1)}=0$  for $t=0,y_n>0$.

The main result of [E] is the following lemma:
\begin{lemma}                           \label{lma:2.1}
Let $T_p,p=1,2,$ be such that the semigeodesic coordinates hold in 
$X_{20}^{(p)},\ \overline{\mathcal{G}}
\subset\Gamma_0$  and $X_{20}^{(p)}$ 
does not
intersect $\partial\Omega^{(p)}\times[0,T_p]$  when  $y_n>0$.  Suppose 
$\hat{\Lambda}^{(1)}=\hat{\Lambda}^{(2)}$ on $\mathcal{G}\times(0,T)$,
where $T=\min(T_1,T_2)$.  Then $u_1^{(1)}=u_2^{(1)}$  for $y'\in \overline{\Gamma},
0\leq y_n\leq \frac{T}{2},y_n\leq t\leq T-y_n$  and
$L_1^{(1)}=L_1^{(2)}$ in $\overline{\Gamma}\times[0,\frac{T}{2}]$.
\end{lemma}

{\bf Remark 2.1}  In [E] we assumed that $\mathcal{D}_p(\mathcal{G}\times(0,T_p))$
does not intersect $\partial\Omega^{(p)}\times[0,T_p]$  when  $y_n>0$.
It is easy to see that it is enough to assume that $X_{20}^{(p)}$  does not 
intersect $\partial\Omega^{(p)}\times[0,T_p]$ for $y_n>0$  since by the domain of
dependence argument the solution 
of $L_1^{(p)}u_p^{(1)}=0,\ u_p^{(1)}=u_{pt}^{(1)}=0$ 
for $t=0,\ u_p^{(1)}=f'$ for $y_n=0,\ u_p^{(1)}=0$ for 
$(\partial\Omega^{(p)}\times[0,T_p])\cap\{y_n>0\}$  restricted
to $X_{20}^{(p)}$ does not change whether $(\partial\Omega^{(p)}\times
[0,T_p])\cap\{y_n>0\}$ intersects  
$\mathcal{D}(\mathcal{G}\times(0,T_p))\setminus X_{20}^{(p)}$ 
or not.
\qed

Denote $D_p=\varphi_p^{-1}(\Gamma\times[0,\frac{T}{2}])\subset\Omega^{(p)}, p=1,2.$
Then $\varphi=\varphi_1^{-1}\circ\varphi_2$  is a diffeomorphism of
$\overline{D}_2$ 
onto  $\overline{D}_1$.  Note that $\hat{g}_1(y',0)=\hat{g}_2(y',0)$  (c.f. 
Remark 2.2 in [E]).  It follows from (\ref{eq:2.4})  that 
  $u_1(\varphi_1^{-1}(y),t)=c_1(y)u_2(\varphi_2^{-1}(y),t)$,
where $c_1(y)\in G_0(\overline{\Gamma}\times[0,\frac{T}{2}])$.  Therefore
\begin{equation}                       \label{eq:2.7}
c_2\circ \varphi\circ  L^{(2)}=L^{(1)}\  \mbox{in}\  \overline{D}_1,
\end{equation}
where  $c_2(x)=c_1(\varphi_1(x))\in G_0(\overline{D}_1)$.  

{\bf Remark 2.2}   
The subset $\Gamma\subset\Gamma_0$ does not have to be small.  
   There may be no global
coordinates $x=(x_1,x_2,...,x_n)$  near $\Gamma$ such that  $x_n=0$ 
is the equation of $\Gamma$.
In this case we take a finite cover $\{V_j\}, j=1,...,N,$
of $\overline{\Gamma}$  and apply Lemma \ref{lma:2.1} to each
$\Gamma_j=\Gamma\cap V_j,\ j=1,...,N$.  Let $T_{pj},
\psi_{pj},\ p=1,2,\ j=1,...,N$,  be the same as $T_p,\varphi_p$  
in Lemma \ref{lma:2.1}.

Let $T_p=\min_{1\leq j\leq N}T_{pj},\ T=\min(T_1,T_2)$.
Maps $\psi_{pj},1\leq j\leq N,$ define  a diffeomorphism  
$\varphi_p$ of $D_p$ onto the manifold $\overline{\Gamma}\times[0,\frac{T}{2}]$
where 
$\overline{D}_p\stackrel{def}{=}\varphi_p^{-1}(\overline{\Gamma}\times[0,\frac{T}{2}])
\subset \overline{\Omega}^{(p)},\ p=1,2,
\ \varphi_p=I$ on $\overline{\Gamma}$.

Note that sets $\mathcal{G},\ X_{20}^{(p)}$ corresponding to 
$\Gamma_j,\ 1\leq j\leq N$  determine sets 
$\mathcal{G},\ X_{20}^{(p)}$  corresponding to the manifold $\Gamma$.
Note also 
that $u_p^{(1)}(y',y_n)$  (see (\ref{eq:2.4}) ) is not a scalar function on
$\Gamma\times[0,T]$,  but a half-density as $(g_p(y',y_n))^{\frac{1}{4}}$. 

Therefore
Lemma \ref{lma:2.1}  implies:
\begin{lemma}                           \label{lma:2.2}
Let $\Gamma$ and $T$ be such that
$\overline{\mathcal{G}}\subset\Gamma_0$,  the semigeodesic coordinates
hold in $X_{20}^{(p)},p=1,2,$ and
$X_{20}^{(p)}\cap(\partial\Omega^{(p)}\setminus
(\Gamma_0\cap\cup_{j=1}^NV_j))=\emptyset$. 
Suppose $\hat{\Lambda}^{(1)}=\hat{\Lambda}^{(2)}$ on $\mathcal{G}
\times (0,T)$.
Denote $\varphi=\varphi_1^{-1}\circ\varphi_2$.  
There exists a gauge transformation $c_2\in  G_0(\overline{D}_1)$ such that
(\ref{eq:2.7}) holds.
or equivalently
\begin{equation}                                 \label{eq:2.8}
L_{1j}^{(1)}=L_{1j}^{(2)}\ \ \mbox{in}\ \ 
\overline{\Gamma}_j\times [0,\frac{T}{2}],\ \ \ 1\leq j\leq N,
\end{equation}
where $L_{1j}^{(p)},p=1,2,$ have the form 
(\ref{eq:2.6}) in local semigeodesic coordinates in 
$\overline{\Gamma}_j\times[0,\frac{T}{2}]$.
For the brevity of notations we shall write 
$L_1^{(1)}=L_1^{(2)}$ in $\overline{\Gamma}
\times[0,\frac{T}{2}]$ instead of (\ref{eq:2.8}).
\end{lemma}
Note that $\overline{D}_p$ is the union of all geodesics in
$\overline{\Omega}^{(p)},p=1,2,$ starting at $\overline{\Gamma}$, 
orthogonal to $\overline{\Gamma}$ and
having the length  $\frac{T}{2}$.

{\bf Remark 2.3}  In \S 4 we shall deal often with the following situation:
Let $V$ be the same neighborhood as in the beginning of this section.  
Suppose $L^{(p)}$ are defined also in $\overline{V},\ p=1,2,$  and 
$L^{(1)}=L^{(2)}$  when $x_n<0$,  i.e. in $V_-\stackrel{def}{=}V\setminus
\overline{\Omega}^{(p)}$.  Let  $y=\varphi_p(x)$  be the same as in
(\ref{eq:2.2}), (\ref{eq:2.3}).  Therefore  
$\varphi=\varphi_1^{-1}\circ\varphi_2$
is the diffeomorphism  of $\overline{D}_2$ onto $\overline{D}_1$.
We shall show that defining 
$\varphi=I$ on $V_-$ we get a diffeomorphism of $\overline{D}_2\cup\overline{V}_-$
onto $\overline{D}_1\cup\overline{V}_-$.

It follows from (\ref{eq:2.2}), (\ref{eq:2.3})  and (\ref{eq:2.7})
that $\varphi(x)$ satisfies the equations:
\[
\|g_1^{jk}\|=\frac{D\varphi}{D x}\ \|g_2^{jk}\|\left(\frac{D\varphi}{D x}\right)^T
\ \ \ \mbox{in}\ \ \ V,
\]
where $\frac{D\varphi}{D x}$
is the Jacobi matrix of $\varphi$.
Since $\varphi=I$  when $x_n=0$  and since $\|g_1^{jk}\|=\|g_2^{jk}\|$
when $x_n\leq 0$  we have that $\varphi=I$  for $x_n\leq 0.$
\qed

Let $B\subset\Gamma\times[0,\frac{T}{2}]$
be a domain such that
$\partial B\cap\Gamma$ 
is open and connected.  Denote  $B_p=\varphi_p^{-1}(B)$.
We assume that $\Omega^{(p)}\setminus \overline{B}_p$ is smooth.
Note that the restriction of $\varphi=
\varphi_1^{-1}\circ\varphi_2$ to   $\overline{B}_2$  maps
$\overline{B}_2$  onto $\overline{B}_1$.   The following extension lemma holds 
(see [Hi], Chapter 8):
\begin{lemma}                            \label{lma:2.3}
There exists a diffeomorphism $\varphi_3$ of $\overline{\Omega}^{(2)}$ onto
$\overline{\Omega}^{(3)}
\stackrel{def}{=}\varphi_3(\overline{\Omega}^{(2)})$
such that
$\varphi_3|_{B_2}=\varphi$  and $\varphi_3|_{\Gamma_0}=I.$
\end{lemma}

{\bf Proof:}  
Since $\overline{D}_2
=\varphi_2^{-1}(\overline{\Gamma}\times[0,\frac{T}{2}])$
and $\varphi=\varphi_1=\varphi_2=I$ on $\overline{\Gamma}$ 
there exists a smooth family $\psi_t,0\leq t\leq 1,$ of embedding of 
$\overline{D}_2$ in $\R^n$  such that $\psi_0=\varphi$
and $\psi_1=I$  on $\overline{D}_2$.
The proof of this fact is similar to the proof of Theorem 8.3.1. in [Hi].
  Denote $A=\overline{B}_2\cup\overline{\Gamma}_0$.
Define $\varphi=I$ on $\overline{\Gamma}_0$.
Denote by $B_2'$ the union of 
a small neighborhood of $\overline{D}_2$
 and a small neighborhood of $\overline{\Gamma}_0$.
We can  extend $\psi_t,0\leq t\leq 1$  from  $\overline{D}_2$  
 to $B_2'$  
preserving 
the properties that $\psi_0=\varphi$ on $A,\ \psi_t=I$  on $\overline{\Gamma}_0,
\  \psi_1=I$ on $B_2',\ \psi_t$ is an embedding of $B_2'$ in
$\R^n$.   Now applying Theorem 8.1.4  in [Hi]
we get that there exists a diffeomorphism  $\varphi_3$ 
of $\R^n$ onto $\R^n,\ \varphi_3=I$ for $|x|>N,\ N$ is large, such that 
$\varphi_3|_A=\varphi$.
Taking the restriction of $\varphi_3$ to $\overline{\Omega}^{(2)}$
we prove Lemma \ref{lma:2.3}.
\qed

Denote by $c_3\in G_0(\overline{\Omega}^{(3)})$ the extension
of   $c_2$ from $\overline{B}_1\cup\overline{\Gamma}_0$
to $\overline{\Omega}^{(3)}$,  where $c_2$ is the same as in 
(\ref{eq:2.7}),  $c_2=1$  on $\Gamma_0\ $,
$\Omega^{(3)}\stackrel{def}{=}\varphi_3(\Omega^{(2)})$.  Let
$L^{(3)}=c_3\circ\varphi_3\circ L^{(2)}$ be 
the differential operator  in $\Omega^{(3)}$.  We have that
$L^{(3)}=L^{(1)}$ in $\overline{B}_1$.
Note that $B_1\subset \Omega^{(1)}\cap\Omega^{(3)}$ and
$\Gamma_0\subset\partial\Omega^{(1)}\cap\partial\Omega^{(3)}$.

Let  $\Gamma_1=(\Gamma_0\setminus(\partial B_1\cap\Gamma_0))\cup
(\partial B_1\setminus\overline{\Gamma}_0)$ 
and let
$\delta=\max_{x\in\overline{B}_1}d_1(x,\Gamma_0)$,   where
$d_1(x,\Gamma_0)$ is the distance in $\overline{B}_1$  from $x\in \overline{B}_1$
to $\Gamma_0$.  Note that $\Gamma_1\subset\partial(\Omega^{(1)}\setminus
\overline{B}_1)\cap\partial(\Omega^{(3)}\setminus\overline{B}_1)$.

Let
$\Lambda^{(1)}$ and $\Lambda^{(3)}$ be the D-to-N operators corresponding 
to $L^{(1)}$  and $L^{(3)}$ in domains $\Omega^{(1)}\times(0,T_0)$ 
and 
$\Omega^{(3)}\times(0,T_0)$,   respectively.  For the simplicity 
of notations we continue to denote by $\Lambda^{(1)},\ \Lambda^{(3)}$
the D-to-N operators corresponding to $L^{(1)}$ and $\L^{(3)}$
in smaller domains $(\Omega^{(1)}\setminus\overline{B}_1)\times
(\delta,T_0-\delta),\ 
(\Omega^{(3)}\setminus\overline{B}_1)
\times(\delta,T_0-\delta)$,  respectively.

The following lemma was proven in [E], Lemma 3.3.
\begin{lemma}                           \label{lma:2.4}
If $\Lambda^{(1)}=\Lambda^{(3)}$ on $\Gamma_0\times(0,T_0)$
then 
the D-to-N operators corresponding  to $L^{(1)},\ L^{(3)}$  on
domains $(\Omega^{(1)}\setminus\overline{B}_1)\times(\delta,T_0-\delta)\ $,
$(\Omega^{(3)}\setminus\overline{B}_1)\times(\delta,T_0-\delta)$  
respectively are equal 
on $\Gamma_1\times(\delta,T_0-\delta)$.
\end{lemma}
Note that the inverse of Lemma \ref{lma:2.4}  is also true:
\begin{lemma}                           \label{lma:2.5}
Suppose 
the D-to-N operators corresponding to $L^{(1)},\ L^{(3)}$
on domains $(\Omega^{(1)}\setminus\overline{B}_1)\times(0,T_0)$,
$(\Omega^{(1)}\setminus\overline{B}_1/\times(0,T_0)$,  respectively,
are equal on  $\Gamma_1\times(0,T_0)$.
Then $\Lambda^{(1)}=\Lambda^{(3)}$ on $\Gamma_0\times(0,T_0)$.
\end{lemma}

{\bf Proof:}  Let
$f$ be a smooth function on $\partial\Omega^{(1)}\times(0,T_0),\ 
\mbox{supp\ }f\subset\Gamma_0\times(0,T_0]$.   Let $u_1$ be the solution 
of the initial-boundary value problem
\begin{eqnarray}
L^{(1)}u_1=0 \ \ \ \mbox{in}\ \ \ \Omega^{(1)}\times(0,T_0),
\ \ \ \ \ \ \ \ \ \ \ \ \ \ \ \ \ \ 
\nonumber
\\
u_1(x,0)=u_{1t}(x,0)=0,\ \ x\in\Omega^{(1)},
\ \ \ \ \ \ u_1|_{\partial\Omega^{(1)}\times(0,T_0)}=f.
\nonumber
\end{eqnarray}
Denote by $f_1$
the restriction of $u_1$ to $\Gamma_1\times(0,T_0)$.
Let $u_3$ be the solution of 
\begin{eqnarray}
L^{(3)}u_3=0 \ \ \ \mbox{in}\ \ \ (\Omega^{(3)}
\setminus\overline{B}_1)\times(0,T_0),
\nonumber
\\
u_3(x,0)=u_{3t}(x,0)=0,\ \ x\in\Omega^{(3)}\setminus\overline{B}_1,
\nonumber
\\
u_3|_{\Gamma_1\times(0,T_0)}=f_1,\ \ u_3=0\ \ \mbox{on}\ \ \
((\partial(\Omega^{(3)}\setminus\overline{B}_1))
\setminus\overline{\Gamma}_1)\times(0,T_0).
\nonumber
\end{eqnarray}
Note that $u_1$ is also zero on $((\partial(\Omega^{(1)}
\setminus\overline{B_1}))\setminus\overline{\Gamma_1})\times(0,T_0)$.
We have that $\Lambda^{(1)}f_1=\Lambda^{(3)}f_1$  on $\Gamma_1\times(0,T_0)$ 
by
the assumption.

Let $\tilde{u}_3=u_3$ in $(\Omega^{(3)}\setminus\overline{B}_1)\times(0,T_0)$
and $\tilde{u}_3=u_1$  in $\overline{B}_1\times(0,T_0)$.   
Since  $\Lambda^{(1)}f_1=\Lambda^{(3)}f_1$
and $L^{(3)}=L^{(1)}$ in $\overline{B_1}$  we get  that $\tilde{u}_3$ is
the solution of $L^{(3)}\tilde{u}_3=0$ in $\Omega^{(3)}\times(0,T_0)$.
Note that $\tilde{u}_3(x,0)=\tilde{u}_{3t}(x,0)=0$ for $x\in \Omega^{(3)}$.  Also
$\tilde{u}_3|_{\Gamma_0\times(0,T_0)}=u_1|_{\Gamma_0\times(0,T_0)}=f,\ 
\ \tilde{u}_3=0$ on
$(\partial\Omega^{(3)}\setminus\Gamma_0)\times
(0,T_0)$.   Moreover,  $\Lambda^{(3)}f=
\Lambda^{(1)}f$  on $\Gamma_0\times(0,T_0)$  since  $u_1=u_3$  in 
$\overline{B_1}\times(0,T_0)$ and $\Lambda^{(1)}=\Lambda^{(3)}$ on 
$\Gamma_1\times(0,T_0)$.
\qed

\section{Refinements of Lemma \ref{lma:2.2}.}
\label{section 3}
\init

In this and the next sections we shall show how to use repeatedly  
Lemmas \ref{lma:2.1} - \ref{lma:2.5} to prove Theorem \ref{theo:1.1}.
We shall prove the following lemma considering,
for the simplicity,
  first the case  when 
$\Gamma=\Gamma_0$ 
and $\partial\Gamma_0=\emptyset$.
\begin{lemma}                            \label{lma:3.1}
Let  $\partial\Gamma_0=\emptyset$  and let 
$T_1$  be such that the semigeodesic coordinates for $L^{(1)}$ hold
in $D_1=\varphi_1^{-1}(\Gamma_0\times[0,\frac{T_1}{2}])$  and
$D_1\cap(\partial\Omega^{(1)}\setminus\Gamma_0)=\emptyset$.  
Let $T_2^*$  be such that      
the semigeodesic coordinates for $L^{(2)}$ hold in 
$\overline{D}_2=\varphi_2^{-1}(\Gamma_0\times[0,\frac{T_2}{2}])$  for any
$T_2<T_2^*$ and $D_2\cap(\partial\Omega^{(2)}\setminus\Gamma_0)=\emptyset$
for any $T_2<T_2^*$.  Suppose that there exists a focal point of a geodesics $\gamma_0$
in $\overline{\Omega}^{(2)}$,  starting at 
$\Gamma_0$, orthogonal to $\Gamma_0$ and such that $y_{n0}=\frac{T_2^*}{2}$,   where
$(y_0',y_{n0})$ are semigeodesic coordinates of the focal point.  Then,  assuming
that $\Lambda^{(1)}=\Lambda^{(2)}$ on $\Gamma_0\times(0,T_0)$ and $T_0>T_1$,
we have that $T_2^*>T_1$ and $L_1^{(1)}=L_1^{(2)}$ in $\Gamma_0\times[0,\frac{T_1}{2}]$.
\end{lemma}
{\bf Proof:}  Consider the bicharacteristics system:
\begin{eqnarray}                         \label{eq:3.1}
\frac{dx_{rj}}{dt}=\frac{\partial H_r(x_r,p_r)}{\partial p_j},\ \ x_{rj}(0)=y_j,
\ 1\leq j\leq n-1,\ r=1,2,
\nonumber
\\
 \frac{dx_{rn}}{dt}=\frac{\partial H_r(x_r,p_r)}{\partial p_n},\ \ x_{rn}(0)=0,
\ \ r=1,2,
\\
\frac{dp_{rj}}{dt}=-\frac{\partial H_r(x_r,p_r)}{\partial x_j},\ \ p_{rj}(0)=0,
\ 1\leq j\leq n-1,\ r=1,2,
\nonumber
\\
\frac{dp_{rn}}{dt}=-\frac{\partial H_r(x_r,p_r)}{\partial x_{n}},\ \ p_n(0)=
\frac{1}{\sqrt{g_r^{nn}(y',0)}},
\ \ r=1,2.
\nonumber
\end{eqnarray}
Here
$H_r(x,p)=\sqrt{\sum_{j,k=1}^n g_r^{jk}(x)p_jp_k}$.
Let $\varphi_r(x)=(\varphi_{r1}(x),...,\varphi_{rn}(x))$ 
be the same as in (\ref{eq:2.2}),  (\ref{eq:2.3}).   Then
$p_{rj}(t)=\frac{\partial \varphi_{rn}(x_r(t))}{\partial x_j},
\ 1\leq j\leq n,
\ r=1,2,$  and (\ref{eq:2.2}) implies  that 
$\frac{d\varphi_{rn}(x(t))}{dt}
=\sum_{j=1}^n p_{rj}(t)\frac{\partial H_r(x_r,p_r)}{\partial p_j}
=H(x_r(t),p_r(t))=1$.

Since $\varphi_{rn}(x(0))=0$  we get that 
\begin{equation}                      \label{eq:3.2}
\varphi_{rn}(x_r(t))=t.
\end{equation}
Note that (\ref{eq:2.3})  implies:
\begin{eqnarray}
\frac{d\varphi_{rj}(x_r(t))}{dt}
=\sum_{k=1}^n\frac{\partial\varphi_{rj}(x_r(t))}{\partial x_k}
\frac{dx_{rk}}{dt}
=\sum_{k=1}^n\frac{\partial\varphi_{rj}(x_r(t))}{\partial x_k}
\frac{\partial H_r(x_r,\frac{\partial \varphi_{rn}}{\partial x})}{\partial p_k}=0,
\nonumber
\\
\ \ \ 1\leq j\leq n-1.
\nonumber
\end{eqnarray}	
Therefore
\begin{equation}                             \label{eq:3.3}
\varphi_{rj}(x_r(t))=\varphi_{rj}(x_r(0))
=\varphi_{rj}(y',0)=y_j,\ 1\leq j\leq n-1,\ r=1,2.
\end{equation}
Denote $y_n=t$. The coordinates $(y_1,...,y_{n-1},y_n)$  are the semigeodesic 
coordinates.
Consider the change of variables
\begin{equation}                          \label{eq:3.4}
x_{j}=x_{rj}(y_1,...,y_n),\ \ \ \ 1\leq j\leq n,\ \ \ r=1,2.
\end{equation}
Assume that the Jacobian
\[
\det\frac{D(x_{r1},...,x_{rn})}{D(y_1,...,y_n)}\neq 0.
\]
It follows from (\ref{eq:3.2}), (\ref{eq:3.3})  that 
(\ref{eq:3.4}) is the inverse to
the change of variables
(\ref{eq:2.1}), since we have :
\begin{eqnarray} 
y_j=\varphi_{rj}(x_r(y_1,...,y_{n-1},t)),
\nonumber
\\
y_n=\varphi_{rn}(x_r(y_1,...,y_{n-1},t))=t.
\nonumber
\end{eqnarray}
Therefore $\frac{D\varphi_r(x)}{Dx}=(\frac{Dx_r}{Dy})^{-1}$.  Note that the 
solution of the system (\ref{eq:3.1}) exists for all $t\in \R$.   Therefore 
$\frac{Dx_r}{Dy}$ is smooth for all $y$.  However the Jacobian  $\det\frac{Dx_r}{Dy}$
may vanish for some $y$.  Such $y$ are called the focal points (caustics)  for 
the Hamiltonian system (\ref{eq:3.1}).   Then $\det\frac{D\varphi_r}{Dx}\rightarrow\infty$
when $y=\varphi_r(x)$  approaches a focal point.  Suppose $T_2^*\leq T_1$.  Then
there exists a focal point  $y^{(0)}=(y_{0}',y_{n0})\in \Gamma_0\times [0,\frac{T_1}{2}]$
for the system (\ref{eq:3.1}) when $r=2$ and $y_{n0}\leq \frac{T_1}{2}$.  
We have near $y^{(0)}$(see (1.8) in [E]):
\begin{equation}                               \label{eq:3.5}
\det\|\hat{g}_r^{jk}(\varphi_r(x))\|=\det\|g_r^{jk}(x)\|
\left(\det\frac{D\varphi_r}{Dx}\right)^2,
\end{equation}
where $\|\hat{g}_r^{jk}(y)\|^{-1}$ is the metric tensor  in semigeodesic
coordinates.  If $y\in\Gamma_0\times[0,\frac{T_2}{2}]$  and $\frac{T_2}{2}$ is close to 
$y_{n0},\ \frac{T_2}{2}<y_{n0}\leq \frac{T_1}{2}$,  we have,  by Lemma \ref{lma:2.2} that
$\hat{g}_1^{jk}=\hat{g}_2^{jk}$.   Consider (\ref{eq:3.5})  when  $r=1$.
Since $y^{(0)}$ is not a focal point for $r=1$ we have that $\det\frac{D\varphi_1}{Dx}$
  is bounded. 
Therefore $\det\|\hat{g}_1^{jk}\|$
is bounded near $y^{(0)}$.  Now consider (\ref{eq:3.5})  for  $r=2$.  Then 
$\left(\det\frac{D\varphi_2}{Dx}\right)^2
=\frac{\det\|\hat{g}_1^{jk}\|}{\det\|g_2^{jk}\|}$  is also bounded when 
$y\rightarrow y^{(0)}$.  Therefore $y^{(0)}$ is not a focal point for 
$r=2$ and this is a contradiction.  Therefore we can take 
$T_2=T_1$ and we have $L_1^{(1)}=L_1^{(2)}$  in 
$\Gamma_0\times[0,\frac{T_1}{2}]$  (c.f. (\ref{eq:2.8})).
\qed

\begin{lemma}                                \label{lma:3.2}
Let 
$\partial\Gamma_0 =\emptyset$ and let
$T_1$ be such that the semigeodesic coordinates for $L^{(1)}$ hold
in $D_1=\varphi_1^{-1}(\Gamma_0\times[0,\frac{T_1}{2}])$  and
$D_1\cap(\partial\Omega^{(1)}\setminus\Gamma_0)=\emptyset$.
Then semigeodesic coordinates for $L^{(2)}$ hold also in
$D_2=\varphi_2^{-1}(\Gamma_0\times[0,\frac{T_1}{2}]),\ 
D_2\cap(\partial\Omega^{(2)}\setminus\Gamma_0)=\emptyset$ and 
$L_1^{(1)}=L_1^{(2)}$ in $\Gamma_0\times[0,\frac{T_1}{2}]$.
\end{lemma}
{\bf Proof:}
Suppose  $D_2$ intersects $\partial\Omega^{(2)}\setminus \Gamma_0$.
  Let  $x^{(0)}$  be  the point in $D_2\cap
(\partial\Omega^{(2)}\setminus\Gamma_0)$ closest to $\Gamma_0$  and
let  $y^{(0)}=(y_0',y_{n0})$ be the semigeodesic coordinates of $x^{(0)}$.
We have by assumption that $y_{n0}\leq\frac{T_1}{2}$.   For any  $\frac{T_2}{2}
<y_{n0}$  we have
by Lemma \ref{lma:3.1}  that  $L_1^{(1)}=L_1^{(2)}$  in $\Gamma_0\times[0,\frac{T_2}{2}]$.
Therefore by the continuity $L_1^{(1)}=L_1^{(2)}$ in $\Gamma_0\times[0,y_{n0}]$.

Denote by $\gamma_0$ the geodesics in $\overline{\Omega}^{(2)}$  starting at 
$x^{(1)}\in \Gamma_0,$
orthogonal to $\Gamma_0$  and reflecting at 
$\partial\Omega^{(2)}\setminus\Gamma_0$ at
point $x^{(0)}$.   Since $x^{(0)}$  is the closest point to $\Gamma_0$
in $D_1$,
the angle of reflection at $x^{(0)}$ is $\frac{\pi}{2}$.  Denote by 
$u_2$  the geometric optics solution depending on large parameter associated 
with $\gamma_0$ and its successive reflections (c.f. [E],  page 824,  
and [E2],  pages 28-29).   
We have  $L^{(2)}u_2=0$ in $\Omega^{(2)}\times(0,T_0),\ u_2=u_{2t}=0$  for
$t=0$,  and $u_2|_{\partial\Omega^{(2)}\times(0,T_0)}=f$  where 
$\mbox{supp\ }f$  is contained in a small neighborhood of $(x^{(1)},t_1)\in 
\Gamma_0\times(0,T_0),t_1>0,$  where $(x^{(1)},t_1)$ is the starting point of 
the broken ray $\gamma_0$. In semigeodesic coordinates $u_2$ is concentrated
(modulo lower order terms)  in a small neighborhood of a broken ray
$y'=y_0',
\ y_n-t=-t_1$  for $t_1<t<y_{n0}+t_1,\ y'=y_0',\ y_n+t=2y_{n0}+t_1$
for $y_{n0}+t_1<t<t_1+2y_{n0},\ y'=y_0',\ y_n-t=-t_1-2y_{n0}$  
for $t>t_1+2y_{n0}$,  etc.

Denote by $u_1$  the solution of $L^{(1)}u_1=0$ in $\Omega^{(1)}\times(0,T_0),\ 
u_1=u_{1t}=0$  for $t=0,\ u_1|_{\partial\Omega^{(1)}\times(0,T_0)}=f$.
Let $u_p^{(1)}$  be related  to $u_p,p=1,2,$  by formulas (\ref{eq:2.4}).
Since $\Lambda^{(1)}=\Lambda^{(2)}$  on $\Gamma_0\times(0,T_0)$  we have that
$u_1^{(1)}$  and $u_2^{(1)}$ have tha same Cauchy data when $y_n=0$.
Since  $L_1^{(1)}=L_1^{(2)}$ in $\Gamma_0\times[0,y_{n0}]$  we get,  by the
unique continuation theorem (see [T]),  that $u_1^{(1)}=u_2^{(1)}$  for 
$y\in\Gamma_0\times[0,y_{n0}],\ 0<t<T_0- y_{n0}$.   However,  $y^{(0)}$
is not a point  of reflection for $u_1^{(1)}$.  Therefore 
$u_1^{(1)}\neq u_2^{(1)}$ 
for $y_{n0}+t_1<t<y_{n0}+t_1+\e,\ y_n<y_{n0}$,  where $\e>0$ is small.
This contradiction proves that $(\partial\Omega^{(2)}\setminus\Gamma_0)\cap
\bar{D}_2=\emptyset$.
\qed

Now we  shall consider the case when $\Gamma\subset\Gamma_0$  
may have a boundary,  i.e. 
when $\partial\Gamma\neq\emptyset$.    
Let $T_1$ be such that all conditions of Lemma \ref{lma:2.2}  are
satisfied.
 Let  $X_{j0}^{(p)},j=1,2,$  be the same as in Lemma \ref{lma:2.2}.
Denote by $\Delta_j^{(p)}$  the projection of $X_{j0}^{(p)}$  on the plane
$t=0,j=1,2,p=1,2.$    Note
that $\Delta_j^{(p)}$ are  contained in the strip $0\leq y_n\leq \frac{T_1}{2}$.
Let $\tilde{\Gamma}\subset\Gamma$  be   such that 
$\overline{\tilde{\Delta}}_{2}^{(1)}\subset \Gamma\times[0,\frac{T_1}{2}]$   where 
$\tilde{\Delta}_{2}^{(p)}$ is the same as 
$\Delta_{2}^{(p)}$
when $\Gamma$ is replaced by $\tilde{\Gamma}$.

We will need the following proposition:
\begin{lemma}                              \label{lma:3.3}
Let $L^{(p)},p=1,2,$  be two operators such that $L_1^{(1)}=L_1^{(2)}$ in
$(\overline{\Gamma}\setminus\tilde{\Gamma})\times[0,\frac{T_1}{2}]$.  Suppose 
$\overline{\tilde{\Delta}}_{2}^{(1)}\subset \Gamma\times[0,\frac{T_1}{2}]$.
Then $\overline{\tilde{\Delta}}_{2}^{(2)}\subset
\overline{\tilde{\Delta}}_{2}^{(1)}$.
\end{lemma}
{\bf Proof:}
Suppose there exists 
$(y_0',y_{n0})\in\overline{\tilde{\Delta}}_{2}^{(2)}$, where
$y_{n0}\leq \frac{T_1}{2},\ (y_0',y_{n0})\notin \overline{\tilde{\Delta}}_{2}^{(1)}$.
Note that $(y_0',y_{n0})\in\overline{\tilde{\Delta}}_{2}^{(2)}$
means that there exists $y_1'\in \overline{\tilde{\mathcal{G}}}$ such
that the shortest path $\gamma$  connecting $(y_1',0)$ and $(y_0',y_{n0})$ has
the length $|\gamma|\leq T_1-y_{n0}$. 
Here $\tilde{\mathcal{G}}$  is the same for $\tilde{\Gamma}$ as
$\mathcal{G}$  for $\Gamma$  (see \S 2).  If $\gamma$ is contained 
completely in $(\overline{\Gamma}\setminus\tilde{\mathcal{G}})\times[0,\frac{T_1}{2}]$
then we must have  $(y_0',y_{n0})\in\overline{\tilde{\Delta}}_{2}^{(1)}$  since
$L_1^{(1)}=L_1^{(2)}$ in $(\overline{\Gamma}\setminus\tilde{\Gamma)}
\times[0,\frac{T_1}{2}]$.
Suppose there is a part of $\gamma$  that does not belong to
$(\overline{\Gamma}\setminus\tilde{\mathcal{G}})\times[0,\frac{T_1}{2}]$.  
Let $(y_2',y_{n2})\in\partial\tilde{\mathcal{G}}\times[0,\frac{T_1}{2}]$
be such that the part $\gamma_1$ of $\gamma$ connecting $(y_2',y_{n2})$ and 
$(y_0',y_{n0})$ is in 
$(\overline{\Gamma}\setminus\tilde{\mathcal{G}})\times[0,\frac{T_1}{2}]$.
Denote by $\gamma_2$  the remaining part of $\gamma$.
Let $\gamma_2'$ be the straight line connecting $(y_{2}',y_{n2})$
and $(y_{2}',0)$.  Since the metric for $L_1^{(2)}$ has  the form
$(dy_n)^2+\sum_{j,k=1}^{n-1}\hat{g}_{2jk}dy_jdy_k$  we have that the length
$|\gamma_2'|$ of $\gamma_2'$ is less or equal than  $|\gamma_2|$. 
Therefore $T-y_{n0}\geq |\gamma|=|\gamma_1|+|\gamma_2|\geq |\gamma_1|+|\gamma_2'|$.
The path $\gamma_2'\cup\gamma_1$ is contained in 
$(\bar{\Gamma}\setminus\tilde{\mathcal{G}})
\times[0,\frac{T_1}{2}]$.   Since  $L_1^{(1)}=L_1^{(2)}$
 in $(\bar{\Gamma}\setminus\tilde{\Gamma})\times[0,\frac{T_1}{2}]$  we have that 
$\gamma_2'\cup\gamma_1$
is contained in $\overline{\tilde{\Delta}}_{2}^{(1)}$.   
Therefore $(y_0',y_{n0})\in \overline{\tilde{\Delta}}_{2}^{(1)}$.  
This contradiction
proves Lemma \ref{lma:3.3}.
\qed

The following lemma generalizes Lemmas \ref{lma:3.1}  and \ref{lma:3.2}
to the case when $\partial\Gamma\neq \emptyset$.
\begin{lemma}                              \label{lma:3.4}
Consider $\Gamma\subset \Gamma_0$ such that $\partial\Gamma\neq \emptyset$.
Let $0<T_1<T_0$  be  such that the semigeodesic coordinates  for $L^{(1)}$  hold in
$\varphi_1^{-1}(\bar{\Gamma}\times[0,\frac{T_1}{2}])$ and 
$(\partial\Omega^{(1)}\setminus\Gamma)\cap \varphi_1^{-1}(\bar{\Gamma}
\times[0,\frac{T_1}{2}])=\emptyset$.  Let $\tilde{\Gamma} \subset\Gamma$  be such that
$\overline{\tilde{\Delta}_{2}^{(1)}}\subset \Gamma\times[0,\frac{T_1}{2}]$.
Suppose that semigeodesic coordinates for $L^{(2)}$ hold in
$(\bar{\Gamma}\setminus\tilde{\Gamma})\times[0,\frac{T_1}{2}],
\ (\partial\Omega^{(2)}\setminus\Gamma)\cap \varphi_2^{-1}((\bar{\Gamma}
\setminus\tilde{\Gamma})
\times[0,\frac{T_1}{2}])=\emptyset$  and
$L_1^{(1)}=L_1^{(2)}$ in $(\bar{\Gamma}\setminus\tilde{\Gamma})
\times[0,\frac{T_1}{2}]$.
Then the semigeodesic coordinates  for $L^{(2)}$ hold in $\bar{\Gamma}
\times[0,\frac{T_1}{2}],\  
(\partial\Omega^{(2)}\setminus\Gamma)\cap \varphi_2^{-1}(\bar{\Gamma}
\times[0,\frac{T_1}{2}])=\emptyset$
and $L_1^{(1)}=L_1^{(2)}$ in $\bar{\Gamma}\times[0,\frac{T_1}{2}]$.
\end{lemma}

{\bf  Proof:}
Since $\overline{\tilde{\Delta}}_{2}^{(2)}\subset
\overline{\tilde{\Delta}}_{2}^{(1)}\subset
\Gamma\times[0,\frac{T_1}{2}]$ (see Lemma \ref{lma:3.3})
 we have that there is no focal points for 
$L_1^{(2)}$ in $\tilde{\Delta}_{2}^{(2)}\setminus(\tilde{\Gamma}\times[0,\frac{T_1}{2}])$.
Repeating the proof of Lemma \ref{lma:3.1} 
with $\Gamma_0$ replaced by $\tilde{\Gamma}$
we get that there is no focal
points for $L_1^{(2)}$  in $\overline{\tilde{\Gamma}}\times[0,\frac{T_1}{2}]$ and
$L_1^{(1)}=L_1^{(2)}$ in $\overline{\tilde{\Gamma}}\times[0,\frac{T_1}{2}]$  
assuming that
$Z\stackrel{def}{=}(\partial\Omega^{(2)}\setminus\Gamma)\cap\varphi_2^{-1}
(\overline{\tilde{\Gamma}}\times[0,\frac{T_1}{2}])=\emptyset$.  
Now we shall show that the set $Z$  
is  empty.   Suppose $Z\neq \emptyset$.
Since  $(\partial\Omega^{(2)}\setminus\Gamma)\cap\varphi_2^{-1}
(\partial\tilde{\Gamma}\times[0,\frac{T_1}{2}])=\emptyset$  
the point $(y_{n0}',y_{n0})\in Z$
closest to $\Gamma$ belongs to $\tilde{\Gamma}\times[0,\frac{T_1}{2}]$, 
 i.e.  $y_0'$ is
an interior point of $\tilde{\Gamma}$.

From this point we can repeat the proof of Lemma \ref{lma:3.2} to get a contradiction.
\qed

\section{The global step.}
\label{section 4}
\init

We start this section with a lemma (Lemma \ref{lma:4.1})
that will play a key role
in the global step of the proof of Theorem \ref{theo:1.1}

Let $O(\delta_1)$ be a ball in $\Gamma_0:O(\delta_1)=\{x'\in\Gamma_0,\ 
d_1(x',x_0')<\delta_1\}$,  where $x_0'\in \Gamma_0,\ d_1(x',x_0')$ is 
the distance on $\Gamma_0$ induced by the metric $\| g_1^{jk}\|^{-1}$.
Let $R$ be the union of geodesics in $\Omega^{(1)}$ with respect to the 
metric $\|g_1^{jk}\|^{-1}$,  starting on $\overline{O}(\delta_1)$,
orthogonal to $\Gamma_0$  and having the lengths $T,\ 2T<T_0$.  We assume that
these geodesics have no focal points in $\overline{R}$ and do not intersect
$\partial\Omega^{(1)}\setminus \overline{O}(\delta_1)$.   Therefore we can
introduce semigeodesic coordinates $y=\psi_1(x)$ in $\overline{R}$ using these
geodesics.  By the continuity the semigeodesic coordinates $y=\psi_1(x)$
hold in a larger domain $R_0\supset R$, i.e. there exists a small $\delta_2>0$
such that $\overline{O}(\delta_1+\delta_2)\subset \Gamma_0$ and $y=\psi_1(x)$
is a diffeomorphism of $\overline{R}_0$  onto 
$\overline{O}(\delta_1+\delta_2)\times[0,T]$.

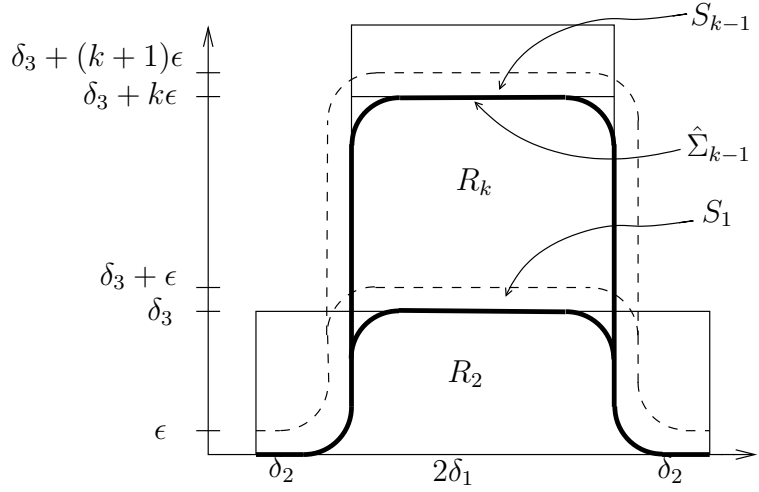
\makefig{Domains $R_k$ and boundaries $\hat{\Sigma}_{k-1}$, $k \ge
  2$. $\hat{\Sigma}_{k-1}$ is drawn in
  bold.}{fig:figure1}{\begin{picture}(0,0)%
\includegraphics{figure1.pstex}%
\end{picture}%
\setlength{\unitlength}{3947sp}%
\begingroup\makeatletter\ifx\SetFigFont\undefined%
\gdef\SetFigFont#1#2#3#4#5{%
  \reset@font\fontsize{#1}{#2pt}%
  \fontfamily{#3}\fontseries{#4}\fontshape{#5}%
  \selectfont}%
\fi\endgroup%
\begin{picture}(5083,3076)(-788,-2232)
\put(1860,-2174){\makebox(0,0)[lb]{\smash{{\SetFigFont{12}{14.4}{\rmdefault}{\mddefault}{\updefault}{\color[rgb]{0,0,0}$2 \delta_1$}%
}}}}
\put(831,-2152){\makebox(0,0)[lb]{\smash{{\SetFigFont{12}{14.4}{\rmdefault}{\mddefault}{\updefault}{\color[rgb]{0,0,0}$\delta_2$}%
}}}}
\put(3272,-2157){\makebox(0,0)[lb]{\smash{{\SetFigFont{12}{14.4}{\rmdefault}{\mddefault}{\updefault}{\color[rgb]{0,0,0}$\delta_2$}%
}}}}
\put(1954,-1553){\makebox(0,0)[lb]{\smash{{\SetFigFont{12}{14.4}{\rmdefault}{\mddefault}{\updefault}{\color[rgb]{0,0,0}$R_2$}%
}}}}
\put(2007,-346){\makebox(0,0)[lb]{\smash{{\SetFigFont{12}{14.4}{\rmdefault}{\mddefault}{\updefault}{\color[rgb]{0,0,0}$R_k$}%
}}}}
\put(3551,-546){\makebox(0,0)[lb]{\smash{{\SetFigFont{12}{14.4}{\rmdefault}{\mddefault}{\updefault}{\color[rgb]{0,0,0}$S_1$}%
}}}}
\put(3480,688){\makebox(0,0)[lb]{\smash{{\SetFigFont{12}{14.4}{\rmdefault}{\mddefault}{\updefault}{\color[rgb]{0,0,0}$S_{k-1}$}%
}}}}
\put(3454,-133){\makebox(0,0)[lb]{\smash{{\SetFigFont{12}{14.4}{\rmdefault}{\mddefault}{\updefault}{\color[rgb]{0,0,0}$\hat{\Sigma}_{k-1}$}%
}}}}
\put(110,-1911){\makebox(0,0)[lb]{\smash{{\SetFigFont{12}{14.4}{\rmdefault}{\mddefault}{\updefault}{\color[rgb]{0,0,0}$\epsilon$}%
}}}}
\put(-230,-941){\makebox(0,0)[lb]{\smash{{\SetFigFont{12}{14.4}{\rmdefault}{\mddefault}{\updefault}{\color[rgb]{0,0,0}$\delta_3+\epsilon$}%
}}}}
\put( 69,-1165){\makebox(0,0)[lb]{\smash{{\SetFigFont{12}{14.4}{\rmdefault}{\mddefault}{\updefault}{\color[rgb]{0,0,0}$\delta_3$}%
}}}}
\put(-788,434){\makebox(0,0)[lb]{\smash{{\SetFigFont{12}{14.4}{\rmdefault}{\mddefault}{\updefault}{\color[rgb]{0,0,0}$\delta_3+(k+1)\epsilon$}%
}}}}
\put(-335,194){\makebox(0,0)[lb]{\smash{{\SetFigFont{12}{14.4}{\rmdefault}{\mddefault}{\updefault}{\color[rgb]{0,0,0}$\delta_3+k\epsilon$}%
}}}}
\end{picture}%
}

Denote by $R_1$  the union of all geodesics in $R_{0}$ 
with the lengths $\delta_3$.   We will choose $\delta_3<\delta_2$. 
   Let  $R_1'$  be the union of geodesics in $\Omega^{(2)}$
with respect to the metric $\|g_2^{jk}\|^{-1}$  starting on 
$\overline{O}(\delta_1+\delta_2)$ 
orthogonal to $\Gamma_0$ and having the lenghts $\delta_3$.
Let  $y=\psi_2(x)$  be  the semigeodesic coordinates on $R'_1$,  i.e.
$\psi_2$  is a diffeomorphism of  $\overline{R}_1'$ 
onto $\overline{O}(\delta_1+\delta_2)
\times[0,\delta_3]$. 

Note that the metrics on $\Gamma_0$  induced by $\|g_1^{jk}\|^{-1}$ and
$\|g_2^{jk}\|^{-1}$  are the same since $\Lambda^{(1)}=\Lambda^{(2)}$ on
$\Gamma_0\times(0,T_0)$  (c.f.,  for example,  Remark 2.2 in [E]).  Note that 
$\psi_3=\psi_1^{-1}\circ\psi_2$ is a diffeomorphism  of $\overline{R}_1'$  onto
$\overline{R}_1$.  Since 
$\delta_3>0$  is small we can apply  Lemma \ref{lma:2.1} or Lemma \ref{lma:2.2}
to get $L^{(1)}=L^{(3)}$   in $\overline{R}_1$  where
$L^{(3)}=c_3\circ \psi_3\circ L^{(2)},\ c_3\in G_0(\overline{R}_1),\ 
\psi_3=I$  and  $c_3=1$ on $O(\delta_1+\delta_2),
 \ L^{(3)}$
is an operator on $\overline{R}_1$.  

Let $\hat{X}_{20}^{(p)}$  be the same
as in  \S 2  when $\Gamma$   is replaced  by $O(\delta_1+\frac{\delta_2}{2})$
and let $\hat{\Delta}_{20}^{(p)}$  be the projection of $\hat{X}_{20}^{(p)}$
on the plane $t=0,p=1,3$.   If $\delta_3$  is much smaller
than $\delta_2$,   we have that 
$\overline{\hat{\Delta}}_{20}^{(1)}\subset O(\delta_1+\delta_2)\times[0,\delta_3]$
and therefore $\overline{\tilde{\Delta}}_{20}^{(3)}\subset
\overline{\tilde{\Delta}}_{20}^{(1)}$ (c.f. Lemma \ref{lma:3.3}).  

Let
$\Sigma_1=\psi_1^{-1}(\sigma_1)$  
where $\sigma_1=\overline{O}(\delta_1+\delta_2)\setminus O(\delta_1)$
when $y_n=0,\ \sigma_1=\partial O(\delta_1)$  when $0\leq y_n\leq \delta_3,\ 
\sigma_1=O(\delta_1)$
 when $y_n=\delta_3$  (see Fig. 1).  Here  $y=\psi_1(x)$ are 
the semigeodesic
coordinates in $\overline{R}_{0}$.  Note that $\sigma_1$  
has edges when $y_n=0,\ y'\in \partial O(\delta_1)$  and when 
$y_n=\delta_3,\ y'\in\partial O(\delta_1)$.
We will smooth $\Sigma_1=\psi_1^{-1}(\sigma_1)$  near these edges to obtain
a smooth surface $\Sigma_1$.  We can arrange the smoothing in such a way that
$\Sigma_1$  and $\hat{\Sigma}_1$ differ only in a small neighborhood of
edges of the size $O(\varepsilon)$ where $0< \e\ll \delta_3$.  Denote by $R_2$
the domain  bounded by $\hat{\Sigma}_1$  and $O(\delta_1+\delta_2)$.  Note
that  $R_2\subset R_1$.  Using Lemma \ref{lma:2.3}  we can extend $\psi_3$
from $\overline{R}_2'\stackrel{def}{=}\psi_3^{-1}(\overline{R}_2)
\subset\overline{\Omega}^{(2)}$ to 
$\overline{\Omega}^{(2)}$  as a diffeomorphism of
$\overline{\Omega}^{(2)}$
onto $\overline{\Omega}^{(3)}\stackrel{def}{=}\psi_3(\Omega^{(2)}),\ \psi_3=I$
on  $\Gamma_0$  and extend $c_3$  from  $\overline{R}_2$  to $\Omega^{(3)}$
as an element of $G_0(\overline{\Omega}^{(3)})$.
It follows  from Lemma  \ref{lma:2.4} that $\Lambda^{(1)}=
\Lambda^{(3)}$  on  $\Sigma_1\times(\delta_3,T_0-\delta_3)$.

Denote by $S_1$  the union of all geodesics in $\overline{\Omega}^{(1)}$
with respect to the metric $\|g_1^{jk}\|^{-1}$,  starting at $\hat{\Sigma}_1$,
orthogonal to $\hat{\Sigma}_1$  and having lengths $\e>0$.  Since $\e>0$
is small (we assume that $0<\e\ll \delta_3\ll \delta_2$)  there is no
focal points in $\overline{S}_1$ and the interior of $S_1$ does not intersect 
$\partial\Omega^{(1)}$.
Since $\e\ll\delta_3$  we can apply Lemma \ref{lma:3.4} to $\overline{S}_1$
and $L^{(1)},L^{(3)}$.  We get that there is  a diffeomorphism $\psi_{41}$
of $\overline{S}_1$ onto 
$\overline{S}_1'\stackrel{def}{=}\psi_{41}(\overline{S}_1)
\subset
\overline{\Omega}^{(3)}$ and  a gauge transformation
$c_{41}\in G_0(\overline{S}_1)$,  such that $\varphi_{41}=I$  on 
$\hat{\Sigma}_1,\ c_{41}=1$  on $\hat{\Sigma}_1$  and  $L^{(1)}=L^{(4)}$
on $S_1$,   where  
$L^{(4)}\stackrel{def}{=}c_{41}\circ\psi_{41}^{-1}\circ L^{(3)}$.

Define $\psi_4'=\psi_{41}^{-1}$  on $\overline{S}_1',\ \psi_4'=I$  on 
$\overline{R}_2,\ c_4'=c_{41}$  on $\overline{S}_1,\ c_4'=1$
on $\overline{R}_2$.   Since  $L^{(1)}=L^{(3)}$  in  $\overline{R}_2$  
we have that
$\psi_4'$  is a $C^\infty$  diffeomorphism on  $\overline{S}_1'\cup\overline{R}_2$
and $c_4\in G_0(\overline{S}_1\cup\overline{R}_2)$ ( c.f. Remark 2.3 ).
Using Lemma \ref{lma:2.3}  we can extend the diffeomorphism 
$\psi_4'$  from  $\overline{S}_1'\cup \overline{R}_2\subset\overline{\Omega}^{(3)}$ 
to $\overline{\Omega}^{(3)}$.   We can also extend $c_4'$ from $\overline{S}_1\cup
\overline{R}_2$  to $\Omega^{(4)}\stackrel{def}{=}\psi_4'(\overline{\Omega}^{(3)})$.

Define  $\psi_4=\psi_4'\circ\psi_3,\ c_4=c_4'(x)c_3(\psi_4^{-1}(x))$.  Then
$L^{(1)}=L^{(4)}$ in $\overline{S}_1\cup\overline{R}_2$  and $L^{(4)}=c_4\circ
\psi_4\circ L^{(2)}$  is an operator in $\Omega^{(4)}$.

Let $\sigma_2$  be the same as $\sigma_1$  with $\delta_3$ replaced 
by $\delta_3+\e$,  let $\Sigma_2=\psi_1^{-1}(\sigma_2)$ and let
$\hat{\Sigma}_2$ be the smoothing of $\Sigma_2$  (see Fig.1).  Denote by
$R_3$  the domain bounded by $\hat{\Sigma}_2$ and $O(\delta_1+\delta_2)$.
Note that  $R_2\subset R_3\subset \overline{S}_1\cup\overline{R}_2$.  
Therefore  $L^{(1)}=L^{(4)}$  
 in 
$\overline{R}_3$.
It follows from Lemma \ref{lma:2.4}  that  $\Lambda^{(1)}=\Lambda^{(4)}$ on 
$\hat{\Sigma}_2\times(\delta_2+\e,T_0-\delta_2-\e)$.

Now repeat the same construction with $\hat{\Sigma}_1$ replaced by $\hat{\Sigma_2}$.

We will get a domain $S_2\subset\overline{\Omega}^{(1)}$  consisting of
all geodesics with respect to the metric $\|g_1^{jk}\|^{-1}$
starting on $\hat{\Sigma}_2$,  orthogonal to $\hat{\Sigma}_2$ and having
the length $\e$,  where  $\e\ll \delta_3$  
is the same as above.
Applying Lemma \ref{lma:3.4}  and Remark 2.3
 we get a diffeomorphism $\psi_5'$ of
$\overline{S}_2'\cup\overline{R}_3$   onto $\overline{S}_2\cup\overline{R}_3$
and a gauge transformation $c_5'\in G_0(\overline{S}_2\cup\overline{R}_3)$
such that $L^{(1)}=L^{(5)}$  
 on
$\overline{S}_2\cup\overline{R}_3$
where
$L^{(5)}=c_5'\circ\psi_5'\circ L^{(4)},\ S_2'\subset\overline{\Omega}^{(4)}$   
is the same as $S_2$ with
respect to the metric $\|g_4^{jk}\|^{-1}$.  Using Lemma \ref{lma:2.3}
we extend
$\psi_5'$  from $\overline{S}_2'\cup\overline{R}_3\subset \overline{\Omega}^{(4)}$
to $\Omega^{(4)}$  as a diffeomorphism of $\overline{\Omega}^{(4)}$  onto
$\overline{\Omega}^{(5)}\stackrel{def}{=}\psi_5'(\overline{\Omega}^{(4)})$.
Define $\psi_5=\psi_5'\circ\psi_4,\ c_5(x)=c_5'(x)c_4(\psi_5^{-1}(x))$.
Then $\psi_5(x)$ is a diffeomorphism of $\overline{\Omega}^{(2)}$ onto
$\overline{\Omega}^{(5)},\ c_5(x)\in G_0(\overline{\Omega}^{(5)},\ 
L^{(5)}=c_5\circ\psi_5\circ L^{(2)}$.

Analogously, let $\sigma_3$  be the same as
$\sigma_1$ with $\delta_2$ replaced by $\delta+2\e,\ \Sigma_3=\psi_1^{-1}(\sigma_3)$ 
and let
$\hat{\Sigma}_3$  be a smoothing of $\Sigma_3$.
Denote by $R_4$ the domain bounded by $\hat{\Sigma}_3$  and $O(\delta_1+\delta_2)$.
Then $\overline{R}_3\subset\overline{R}_4\subset(\overline{S}_2\cup\overline{R}_3)$
and $L^{(1)}=c_6\circ\psi_6\circ L^{(2)}$  on $\overline{R}_4$,  where $\psi_6$
is a diffeomorphism of $\overline{\Omega}^{(2)}$ onto $\overline{\Omega}^{(6)}
\stackrel{def}{=}\psi_6(\overline{\Omega}^{(2)})$  and 
$c_6\in G_0(\overline{\Omega}^{(6)})$.
After $\frac{T-\delta_3}{\e}=N$  steps we get  a domain  $\overline{R}_N
\supset\overline{R}$,  a diffeomorphism $\psi_{N+2}$  of $\overline{\Omega}^{(2)}$
onto $\overline{\Omega}^{(N+2)}$
and a gauge transformation $c_{N+2}\in G_0(\overline{\Omega}^{(N+2)})$  such
that $L^{(1)}=c_{N+2}\circ\psi_{N+2}\circ L^{(2)}$  
 on $\overline{R}_N,\ \psi_{N+2}=I$  
and $c_{N+2}=1$  on $\Gamma_0$.

We proved the following lemma:
\begin{lemma}                             \label{lma:4.1}
Let 
$O'\subset \Gamma_0$  be a small neighborhood of $O(\delta_1+\delta_2)$.
  Suppose $\Lambda^{(1)}=\Lambda^{(2)}$  
on $O'\times(0,T_0)$.  Let  $R\subset\Omega^{(1)}$  be the same as above and
$T<\frac{T_0}{2}$.   Then  there exists a diffeomorphism $\varphi$ of $R$  onto
$\overline{R'}\stackrel{def}{=}\varphi(\overline{R})\subset\overline{\Omega}^{(2)}$
and a gauge tranformation $c\in G_0(\overline{R})$  such that $L^{(1)}=
c\circ\varphi^{-1}\circ L^{(2)}$    
  on $\overline{R},\ \varphi=I$  and 
$c=1$ on $O(\delta_1)$.
\end{lemma}

{\bf Remark 4.1}  In order to prove Lemma \ref{lma:4.1}   we used 
$L^{(p)},p=1,2$  only in a small neighborhoods of $\overline{R}$ and
$\overline{R'}$.  The properties of $L^{(p)},p=1,2,$ 
outside of these neighborhoods play no role.  For example,  $L^{(p)},p=1,2,$
are not required to be self-adjoint outside  of neighborhoods  of $R$ and
$R'$.

{\bf Remark 4.2}  The following generalization of Lemma \ref{lma:4.1}
holds:  

Let $x^{(0)}\in\Gamma_0$ and let $x^{(1)}$  be an arbitrary point in $\Omega^{(1)}$.   
Let  $\gamma$  be an arbitrary curve in $\Omega^{(1)}$  connecting points $x^{(0)}$
and $x^{(1)}$  and having the length less than $\frac{T_0}{2}$.  There exists
a neighborhood $V_0\subset \Omega^{(1)}$   of $\gamma$,  a diffeomofphism
$\psi_0$ of $\overline{V}_0$ onto $\overline{V}_0'\stackrel{def}{=}
\psi_0(\overline{V}_0)\subset\Omega^{(2)}$ and
there exists a gauge $c_0(x)\in G_0(\overline{V}_0)$  such that 
$L^{(1)}=c_0\circ\psi_0^{-1}\circ L^{(2)}$  
in $\overline{V}_0,\ \psi_0=I$  and $c_0=1$ on $\overline{V}_0\cap\Gamma_0$
assuming that $\Lambda^{(1)}=\Lambda^{(2)}$  on $\tilde{O}\times(0,T_0)$  
where $\tilde{O}\subset\Gamma_0$  is a small neighborhood 
of $\overline{V}_0\cap\Gamma_0$.

To prove this result we approximate $\gamma$  by a piece-wise smooth curve 
consisting of geodesic segments and prove Lemma \ref{lma:4.1} successively 
for each geodesic segment (see similsr arguments below).
\qed

Now we actually start the global construction (c.f. [KKL] and [KKL1]).  Since
$\overline{\Omega}^{(1)}$ is compact there exists $\delta_0>0$  such that for any 
point $x^{(0)}\in \overline{\Omega}^{(1)}$ the geodesics starting at $x^{(0)}$ 
form a local
system of coordinates in $B(x^{(0)},\delta_0)\setminus B(x^{(0)},\delta)$ for
any $0<\delta<\delta_0$.  Here $B(x^{(0)},r)=\{x:d(x,x^{(0)})<r\}$ is
a ball of radius $r$.  We assume that the metric $\|g_1^{jk}\|^{-1}$ is extended
to a $\delta_0$-neighborhood of $\overline{\Omega}^{(1)}$.
Therefore $B(x^{(0)},\delta_0)$ make sense when $B(x^{(0)},\delta_0)\cap\partial\Omega^{(1)}
\neq\emptyset$.  Let $\e_1>0$ be such that the semigeodesic coordinates with
respect to $\partial\Omega^{(1)}$ hold in $2\e_1$-neighborhood of $\partial\Omega^{(1)}$ 
and the interior of this neighborhood does not 
intersect $\partial\Omega^{(1)}$.
For each $x^{(0)}\in\partial\Omega^{(1)}$ consider a ball $B(x^{(0)},\e_1)$.
For each $x^{(0)}\in\Omega^{(1)}$ consuder a ball $B(x^{(0)},r)$
such that $r<\frac{\delta_0}{2}$ and $\overline{B}(x^{(0)},r)\cap\partial\Omega^{(1)}
=\emptyset$.
Such balls form an open cover of $\overline{\Omega}^{(1)}$  and since
$\overline{\Omega}^{(1)}$ is compact there exists a finite subcover
$\{B(x^{(k)},r_k)\},k=1,...,N$.
Denote by $\Omega_{\e_1}$ the union of all balls $B(x^{(k)},r_k)$ such that
$\overline{B}(x^{(k)},r_k)\cap\partial\Omega^{(1)}=\emptyset$.

Let $\Gamma$ and $T$ be the same as in Lemma \ref{lma:2.2}.  Repeating the proof 
of Lemmas \ref{lma:2.2} and \ref{lma:2.3}
and using th same notations we get a domain $B_1\subset\Omega^{(1)},\ 
\partial B_1\cap \Gamma\neq\emptyset$ and connected, a diffeomorphism $\varphi_3$ of
$\overline{\Omega}^{(2)}$ onto $\overline{\Omega}^{(3)}\stackrel{def}{=}
\varphi_3(\overline{\Omega}^{(2)})$ and a gauge 
transformation $c_3\in G_0(\overline{\Omega}^{(3)})$
such that $L^{(3)}=L^{(1)}$ in $B_1\subset \Omega^{(1)}\cap\Omega^{(3)}$  
where $L^{(3)}=c_3\circ\varphi_3\circ L^{(2)}$  is an operator in 
$\Omega^{(3)},\ \varphi_3=I$ on $\Gamma_0$.
Since $\varphi_3=I$ on $\Gamma_0$  we have that $\Lambda^{(3)}=\Lambda^{(1)}$
on $\Gamma_0\times(0,T_0)$  where $\Lambda^{(3)}$ is the D-to-N for $L^{(3)}$.
Let $\Gamma_1=
(\Gamma_0\setminus(\Gamma_0\cap\partial B_1))\cup(\partial B_1\setminus \Gamma_0)$.
Note that $\Gamma_1\subset\partial(\Omega^{(1)}\setminus \overline{B}_1)\cap
\partial(\Omega^{(3)}\setminus \overline{B}_1)$.

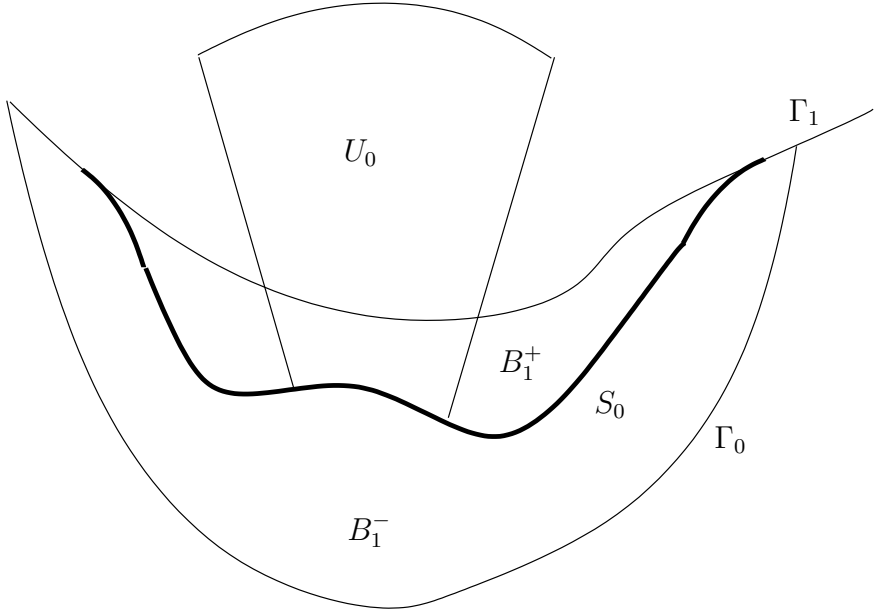
\makefig{The domain $U_0$. The boundary $S_0$ is drawn bold.}{fig:figure2}{\begin{picture}(0,0)%
\includegraphics{figure2.pstex}%
\end{picture}%
\setlength{\unitlength}{3947sp}%
\begingroup\makeatletter\ifx\SetFigFont\undefined%
\gdef\SetFigFont#1#2#3#4#5{%
  \reset@font\fontsize{#1}{#2pt}%
  \fontfamily{#3}\fontseries{#4}\fontshape{#5}%
  \selectfont}%
\fi\endgroup%
\begin{picture}(5464,3827)(-760,-2508)
\put(1392,-2066){\makebox(0,0)[lb]{\smash{{\SetFigFont{12}{14.4}{\rmdefault}{\mddefault}{\updefault}{\color[rgb]{0,0,0}$B_1^-$}%
}}}}
\put(2345,-1013){\makebox(0,0)[lb]{\smash{{\SetFigFont{12}{14.4}{\rmdefault}{\mddefault}{\updefault}{\color[rgb]{0,0,0}$B_1^+$}%
}}}}
\put(2945,-1273){\makebox(0,0)[lb]{\smash{{\SetFigFont{12}{14.4}{\rmdefault}{\mddefault}{\updefault}{\color[rgb]{0,0,0}$S_0$}%
}}}}
\put(3698,-1493){\makebox(0,0)[lb]{\smash{{\SetFigFont{12}{14.4}{\rmdefault}{\mddefault}{\updefault}{\color[rgb]{0,0,0}$\Gamma_0$}%
}}}}
\put(4165,567){\makebox(0,0)[lb]{\smash{{\SetFigFont{12}{14.4}{\rmdefault}{\mddefault}{\updefault}{\color[rgb]{0,0,0}$\Gamma_1$}%
}}}}
\put(1371,314){\makebox(0,0)[lb]{\smash{{\SetFigFont{12}{14.4}{\rmdefault}{\mddefault}{\updefault}{\color[rgb]{0,0,0}$U_0$}%
}}}}
\end{picture}%
}

Take arbitrary point $x^{(0)}\in B_1$ such that $d(x^{(0)},\Gamma_1)=2\e_2$ 
where $\e_2$ is much smaller than $\delta_0$.  We have that  $B(x^{(0)},\e_2)
\subset B_1$.  Pick some geodesics $\gamma_0$ starting at $x^{(0)}$ and denote by
$U_0$ the part of $B(x^{(0)},\delta)\setminus\overline{B}(x^{(0)},\e_2)$ 
consisting of  all geodesics $\gamma$  starting at $x^{(0)}$ and having an angle 
less than $\pi-\e_3$ with $\gamma_0$ at $x^{(0)}$.
We choose $\delta\leq\delta_0,\gamma_0$ and $\e_3$ such that $\overline{U}_0
\cap\partial\Omega^{(1)}=\emptyset$.
Denote by $S_0$ a smooth surface in $B_1$  that contains $\overline{U}_0
\cap\partial B(x^{(0)},\e_2)$  and divides $B_1$ in two domains,  $B_1^+$  and $B_1^-$,
where $B_1^+$  is bounded by $\Gamma_1$ and $S_0$,  and
$B_1^-=B_1\setminus\overline{B}_1^+$  (see Fig.2).  We assume that 
$\Omega^{(1)}\setminus \overline{B}_1^-$  has a smooth boundary that includes $S_0$.
Denote $\Gamma_2=S_0\cup (\Gamma_1\setminus\partial B_1^+)$.  
Since $L^{(3)}=L^{(1)}$  in $\overline{B}_1^-$  and $\Lambda^{(3)}=\Lambda^{(1)}$
on $\Gamma_0\times(0,T) $  we get from the Lemma \ref{lma:2.4}  that
$\Lambda^{(3)}=\Lambda^{(1)}$ on $\Gamma_2\times(\delta',T-\delta')$,
where $\delta'=\max_{x\in\bar{B}_1^-}d_1(x,\Gamma_0)$  and we consider 
the D-to-N operators corresponding to $L^{(1)},L^{(3)}$ in domains 
$(\Omega^{(1)}\setminus \overline{B}_1^-)\times(\delta',T_0-\delta'),\ 
(\Omega^{(3)}\setminus \overline{B}_1^-)\times(\delta',T_0-\delta')$.

Now apply Lemma \ref{lma:4.1}  to $U_0\subset(\Omega^{(1)}\setminus B_1^-)$ 
instead of $R$.
Let  $U_0'$ be the union of all geodesics in $\Omega^{(3)}$ with respect 
to the metric $\|g_3^{jk}\|^{-1}$  starting at $\overline{U}_0
\cap\partial B(x^{(0)},\e_2)$,  orthogonal to $\partial B(x^{(0)},\e_2)$ and
having the length $\delta-\e_2$.  Note that $B(x^{(0)},\e_2)\subset B_1\subset
\Omega^{(1)}\cap\Omega^{(3)}$.
Since $\Lambda^{(1)}=\Lambda^{(3)}$ on $\Gamma_2\times(\delta',T_0-\delta')$,
the Lemma \ref{lma:4.1}
implies that there exists a diffeomorphism $\varphi$  of $\overline{U}_0'
\subset\overline{\Omega}^{(2)}$ onto
$\overline{U}_0$  and there exists $c(x)\in C^\infty(\overline{U}_0),\ c(x)\neq 0$ in 
$\overline{U}_0$  such that $c\circ\varphi\circ L^{(3)}=L^{(1)}$ in 
$\overline{U}_0,\ \varphi=I$  and  $c=1$  on $S_0\cap\overline{U}_0$.
Define $\varphi=I$  on $\overline{B}_1\setminus U_0$ and $\Gamma_0$.
Also define $c=1$  on $\overline{B}_1\setminus U_0$  and $\Gamma_0$.
Since $L^{(3)}=L^{(1)}$  in $\overline{U}_0\cap\overline{B}_1$ we have that
$\varphi=I$ and $c=1$ in $\overline{U}_0\cap\overline{B}_1$.   Therefore $\varphi$ and
$c$ are $C^\infty$ on $\overline{U}_0\cup\overline{B}_1\cup\Gamma_0$
(c.f. Remark 2.3).
Applying Lemma \ref{lma:2.3}  we can extend $\varphi$ from
$\overline{U}_0\cup\overline{B}_1\cup\Gamma_0$
to $\overline{\Omega}^{(3)}$ as a diffeomorphism of $\overline{\Omega}^{(3)}$
onto $\overline{\Omega}^{(4)}\stackrel{def}{=}\varphi(\overline{\Omega}^{(3)})$
and extend $c$ from
$\overline{U}_0\cup\overline{B}_1\cup\Gamma_0$
to $\overline{\Omega}^{(4)}$ as an element of $G_0(\overline{\Omega}^{(4)})$.

Let $\varphi_4=\varphi\circ\varphi_3$.  Then $\varphi_4$ is a diffeomorphism of
$\overline{\Omega}^{(2)}$ onto $\overline{\Omega}^{(4)}$.  Let $c_4(x)=
c(x)c_3(\varphi_4^{-1}(x))\in G_0(\Omega^{(4)}$.
Then  we have $L^{(4)}=c_4\circ\varphi_4\circ L^{(2)}$  is an operator in
$\Omega^{(4)}$  equal to $L^{(1)}$ in $\overline{U}_0\cup\overline{B}_1 
 \subset\overline{\Omega}^{(1)}\cap\overline{\Omega}^{(4)}$.

Therefore applying Lemma \ref{lma:4.1}  to $U_0$  we gained that
$\overline{B}_1$ is replaced by 
a  larger domain
$\overline{U}_0\cup\overline{B}_1$.

Taking  a point $x^{(1)}\in U_0\cup B_1$  instead of $x^{(0)}$  we can construct 
a domain $U_1$  similar to $U_0$.
For the brevity we shall call by $U$-type domains the domains similar to $U_0$.
We shall show that adding  a finite number of $U$-type domains we can cover 
$\overline{\Omega}^{(1)}$. 

Take any ball $B(x^{(p)},r_p)\subset \Omega_{\e_1}$.
Since $\overline{\Omega}^{(1)}$  is connected, the point $x^{(p)}$ 
can be connected with $x^{(0)}\in B_1$ by a broken
geodesics (c.f. [KKL]),  more exactly,  there exist points $y_1=x^{(0)},
y_2,...,y_{N_1}=x^{(p)}$  such that $y_k$  and $y_{k+1},1\leq k\leq N_1-1$,
can be connected by a geodesics of length  $\leq \frac{\delta_0}{2}$.  Using
a sequence of $U$-type domains we can cover this broken geodesics including
$x^{(p)}$.  Adding more $U$-type domains if needed we can cover 
$B(x^{(p)},r_p)$ too.  We can do this with any ball 
$B(x^{(p)},r_p)\subset \Omega_{\e_1}$.  Therefore inserting $M$ 
of $U$-type domains
$U_1,...,U_M$  in $\Omega^{(1)}$  we get a sequences of domains 
$\Omega^{(k)},4\leq k\leq M+4$, diffeomorphisms $\varphi_k$  of 
$\overline{\Omega}^{(2)}$  onto $\overline{\Omega}^{(k)}$
and gauge transformations $c_k\in G_0(\overline{\Omega}^{(k)})$
such that $L^{1)}=L^{(k)}$  
in $\overline{B}_1\cup(\cup_{j=0}^{k-4}\overline{U}_j),\ 4\leq k\leq M+4$,  where 
$L^{(k)}=c_k\circ\varphi_k\circ L^{(2)}$  are operators in $\Omega^{(k)},\ 
B_1\cup(\cup_{j=0}^{k-4}U_j)\subset\Omega^{(1)}\cap\Omega^{(k)},\ 
\Lambda^{(1)}=\Lambda^{(k)}$  on $\Gamma_0\times(0,T_0)$   
$\ 4\leq k\leq M+4$.

Let $M$ be such that $\Omega_{\e_1}\subset B_1\cup(\cup_{j=0}^{M}U_j)\subset\Omega^{(1)}\cap\Omega^{(M+4)}$.
We have 
$L^{(1)}=L^{(M+4)}$
 in
$\overline{\Omega}_{\e_1}$,  
 and $\Lambda^{(1)}=\Lambda^{(M+4)}$ on $\Gamma_0\times
(0,T_0)$. 

Note that we choose $U_j,0\leq j\leq M$, taking into account the geometry of 
$\Omega^{(1)}$ and the metric $\|g_1^{jk}\|^{-1}$, 
regardless of the geometry of $\Omega^{(p)}$  and the metric $\|g_p^{jk}\|^{-1}$,
$\ 2\leq p$.

Now consider the cover of $\overline{\Omega}^{(1)}\setminus \Omega_{\e_1}$  by
$U$-type domains.  Let $\alpha=\partial\Omega^{(1)}\cap B(x^{(k)},\e_1)$,  
where $\{B(x^{(k)},\e_1)\}$  is a finite cover of $\partial\Omega^{(1)}$.
  Denote by
$U_{M+1}$  the union of all geodesics with respect to the metric $\|g_1^{jk}\|^{-1}$
starting on   $\overline{\alpha}$, orthogonal to 
$\overline{\alpha}$  and having the same lengths $2\e_1$.
Let $W_{M+1}$ be the set of endpoints of these geodesics.  It follows from
the definition of $\Omega_{\e_1}$  that $W_{M+1}$  is located inside 
$\Omega_{\e_1}$.

Denote by $U_{M+1}'$ the union of all geodesics with respect to metric
$\|g_{M+4}^{jk}\|^{-1}$ corresponding to $L^{(M+4)}$  that start on $W_{M+1}$,
orthogonal to $W_{M+1}$ and having the lengths $2\e_1$.  Note that
$W_{M+1}\subset\Omega_{\e_1}\subset\Omega^{(M+4)}$. 

Let $\tilde{W}_{M+1}$  bea surface in $\Omega_{\e_1}$  containing 
$W_{M+1}$  and  similar to $S_0$ in the case of domain $U_0$.  Since 
$\Lambda^{(1)}=\Lambda^{(M+4)}$  on $\Gamma_0\times(0,T_0)$  and
$L^{(1)}=L^{(M+4)}$  on $\overline{\Omega}_{\e_1}$,  Lemma \ref{lma:2.4}
implies that  $\Lambda^{(1)}=\Lambda^{(M+4)}$
on $\tilde{W}_{M+4}\times(T',T_0-T')$,  
where $T'=\max_{x\in\overline{\Omega}_{\e_1}} d_1(x,\Gamma_0)$.
Now we can repeat for $U_{M+1}$  the same arguments as for $U_0$.
Applying Lemma \ref{lma:4.1} we get
that there exists 
a diffeomorphism $\varphi'_{M+5}$ of $\overline{U}_{M+1}'$ onto $\overline{U}_{M+1},\ 
\varphi_{M+5}'=I$  on $W_{M+1}$ and $c_{M+5}'(x)\in G_0(\overline{U}_{M+1}),\ 
c'=1$ on $W_{M+1}$ such that 
\begin{equation}                           \label{eq:4.1}
c_{M+5}'\circ\varphi_{M+5}'\circ L^{(M+4)}=L^{(1)}
 \ \mbox{in\ } U_{M+1}.
\end{equation}
  Since $L^{(M+4)}=L^{(1)}$ in $\Omega_{\e_1}$ we have that 
$\varphi_{M+5}'=I$ on $\overline{U}_{M+5}\cap\overline{\Omega}_{\e_1}$.
Therefore $\varphi_{M+5}'$  is a diffeomorphism on 
$\overline{U}_{M+5}\cup\overline{\Omega}_{\e_1}$ after extending $\varphi_{M+5}'$ as 
$I$ on $\overline{\Omega}_{\e_1}$ (c.f. Remark 2.3).  
Analoguously taking $c_{M+5}'=1$ on 
$\overline{\Omega}_{\e_1}$ we get a $C^\infty$-function on 
$\overline{\Omega}_{\e_1}\cup\overline{U}_{M+1}$.
Note that Lemma \ref{lma:4.1}  implies that $U_{M+1}'\subset \Omega^{(M+4)}$.

We shall consider first the case when $\alpha\cap\Gamma_0\neq\emptyset$.
We shall prove that $\varphi_{M+5}'=I,\ c_{M+5}'=1$  on $\alpha_{1}\stackrel{def}{=}
\alpha\cap\Gamma_0$.

Let $g$ be arbitrary smooth function with the support in $\overline{\alpha}_{1}
\times(T',T_0-T')$.

Let $u_1$ be the solution of $L_1u_1=0$ on $\Omega^{(1)}\times(0,T_0),\ 
u_1=u_{1t}=0$  for $t=0,\ u_1|_{\Gamma_0\times(0,T_0)}=g$.  Also let $u_{M+4}(x)$
be the solution of $L^{(M+4)}u_{M+4}=0$  in $\Omega^{(M+4)}\times(0,T_0),\ 
u_{M+4}=\frac{\partial u_{M+4}}{\partial t}=0$  for $t=0$  and
$u_{M+4}|_{\Gamma_{0}\times(0,T_0)}=g$.
Since $\Lambda^{(1)}=\Lambda^{(M+4)}$ on $\Gamma_0\times(0,T)$ and since
$L^{(M+4)}=L^{(1)}$ in $\overline{\Omega}_{\e_1}$ we get, by the unique continuation
theorem (c.f [T]), that $u_1=u_{M+4}$ in $\overline{\Omega}_{\e_1}\times
(T',T_0-T')$.

It follows from (\ref{eq:4.1}) that $u_1(x,t)$ and 
$u_{M+5}(x,t)=c_{M+5}'(x)u_{M+4}((\varphi_{M+5}')^{-1}(x),t)$
satisfy the same equation $L^{(1)}u=0$ in $U_{M+1}\times(T',T_0-T'))$.
Since $u_1(x)=u_{M+4}(x),\ \varphi_{M+5}'=I$  and $c_{M+5}'=1$  on 
$\Omega_{\e_1}\times(T',T_0-T')$ 
we have that $u_1$  and $u_{M+5}$  have the same Cauchy data
on $\tilde{W}_{M+1}\times(T',T_0-T')$.

We
 have,  using again the unique continuation theorem that
\begin{equation}                                \label{eq:4.2}
c_{M+5}'(x)u_{M+4}((\varphi_{M+5}')^{-1}(x),t)=u_1(x,t)
\  \mbox{in  \ }\overline{U}_{M+1}\times(T_1+2\e_1,T_0-T_1-2\e_1).
\end{equation}
Taking the restriction of (\ref{eq:4.2})  to 
$\overline{\alpha}_1\times(T'+2\e_1,T_0-T'-2\e_1)$ we obtain
\begin{equation}                                \label{eq:4.3}
c_{M+5}'(x)g((\varphi_{M+5}')^{-1}(x),t)=g(x,t)  \ \mbox{in\  }
\overline{\alpha}_{1}\times(T_1+2\e_1,T_0-T_1-2\e_1).
\end{equation}

Since $g$ is arbitrary we get that $\varphi_{M+5}'=I$ on $\overline{\alpha}_{1},\ 
c_{M+5}'(x)=1$ on $\overline{\alpha}_{1}$.  Therefore,   we take $\varphi_{M+5}'=I,\ 
c_{M+1}'(x)=1$ on $\Gamma_0$  and get a $C^\infty$-functions 
on $\Gamma_0\cup\overline{\Omega}_{\e_1}\cup\overline{U}_{M+1}$.
There is no any obstruction to the extension $\varphi_{M+5}'=I,\ 
c_{M+5}'=1$ on $\Gamma_0$  in the case  when
$\alpha\cap\Gamma_0=\emptyset$.  Thus in both cases
applying Lemma \ref{lma:2.3}  we get a diffeomorphism 
$\varphi_{M+5}'$ of $\Omega^{(M+4)}$ onto 
$\overline{\Omega}^{(M+5)}\stackrel{def}{=}
\varphi_{M+5}'(\overline{\Omega}^{(M+4)})$ and
$c_{M+5}'\in G_0(\Omega^{(M+5)})$ such that $\varphi_{M+5}'=I$ on $\Gamma_0$.

As in the case of the $U$-type domain $U_0$, denote $\varphi_{M+5}=
\varphi_{M+5}'\circ\varphi_{M+4},\ c_{M+5}=c_{M+5}'c_{M+4}(\varphi_{M+4}^{-1}(x))$.
Then we get that 
\[
L^{(M+5)}=L^{(1)}\ \ \ \mbox{in\ \ \ } \overline{\Omega}_{\e_1}\cup\overline{U}_{M+1}
\subset\overline{\Omega}^{(1)}\cap\overline{\Omega}^{(M+5)},
\]
where $L^{(M+5)}=c_{M+5}\circ\varphi_{M+5}\circ L^{(2)}$
is an operator on $\overline{\Omega}^{(M+5)}$,
$\ \varphi_{M+5}$ is
a diffeomorphism of $\overline{\Omega}^{(2)}$ onto 
$\overline{\Omega}^{(M+5)}$, $\varphi_{M+5}=I$  on $\Gamma_0$,  $c_{M+5}\in
G_0(\overline{\Omega}^{(M+5)})$.
In particular,  we have $\alpha\subset\partial\Omega^{(1)}\cap\partial\Omega^{(M+5)}$
and $\varphi_{M+5}^{-1}$  maps $\alpha\subset\partial\Omega^{(1)}$  onto
$\alpha'\in\partial\Omega^{(2)}$.
Repeating the same construction with each of $B(x^{(k)},\e_1)$ such that
$B(x^{(k)},\e_1)\cap\partial\Omega^{(1)}\neq\emptyset$
 we get $\overline{\Omega}_{\e_1}\cup
(\cup_{j=1}^{M_1}\overline{U}_{M+j})=\overline{\Omega}^{(1)},\ \
\overline{\Omega}_{\e_1}\cup
(\cup_{j=1}^{M_1}\overline{U}_{M+j})\subset\overline{\Omega}^{(M+M_1+4)}$
and $L^{(1)}=c_{M+M_1+4}\circ\varphi_{M+M_1+4}\circ L^{(2)}$  in
$\overline{\Omega}_{\e_1}\cup
(\cup_{j=1}^{M_1}\overline{U}_j)$,
where 
$\overline{\Omega}^{M+M_1+4)}\stackrel{def}{=}
\varphi_{M+M_1+4}(\overline{\Omega}^{(2)}),\
\varphi_{M+M_1+4}(x)$  is a diffeomorphism,  $\varphi_{M+M_1+4}=I$
on $\Gamma_0,\ c_{M+M_1+4}\in G_0(\Omega^{M+M_1+4)})$.

We claim that $\Omega^{(M+M_1+4)}=\Omega^{(1)}$. If 
$\overline{\Omega}^{(M+M_1+4)}  \neq \overline{\Omega}^{(1)}$
then there exists an interior point of $\Omega^{(M+M_1+4)}$
that is a boundary point of $\Omega^{(1)}$.
It follows from the arguments similar to the proof of Lemma \ref{lma:3.2}
that this is impossible.
\qed

\end{document}